\documentclass[11pt,letterpaper]{article} 

\usepackage{cite}

\usepackage{amsthm,amsmath,amssymb,amsfonts}
\usepackage{bm}  % bold greek letters

\usepackage{algorithm}
\usepackage{algpseudocode}
\usepackage{url}

\usepackage{fullpage}

\usepackage{graphicx}
\usepackage{xcolor}

\usepackage{array}
\usepackage{multirow}

\makeatletter
\let\MYcaption\@makecaption
\makeatother

\usepackage[font=footnotesize]{subcaption}

\makeatletter
\let\@makecaption\MYcaption
\makeatother

\theoremstyle{plain}
\newtheorem{theorem}{Theorem}%[section]
\newtheorem{lemma}{Lemma}%[section]

\theoremstyle{definition}

\newtheorem{definition}{Definition}
\newtheorem{assumption}{Assumption}
\newtheorem{remark}{Remark}

\def\({\left(}
\def\){\right)}
\def\[{\left[}
\def\]{\right]}

\def\sec#1{Section~\ref{sec:#1}}

\title{\LARGE \bf Explicit Distributed and Localized Model Predictive Control via System Level Synthesis}

\author{Carmen Amo Alonso\thanks{C. Amo Alonso is with the  Computing and Mathematical Sciences Department at California Institute of Technology, Pasadena, CA.
        {\tt\small camoalon@caltech.edu}}, Nikolai Matni\thanks{N. Matni is  with the Department of Electrical and Systems Engineering at the University of Pennsylvania, Philadelphia, PA.    {\tt\small nmatni@seas.upenn.edu}}, and James Anderson% <-this % stops a space
\thanks{J. Anderson is with the Department of Electrical Engineering and the Data Science Institute at Columbia University, New York, NY. {\tt\small james.andersosn@columbia.edu}}}

\begin{document}

\maketitle
\thispagestyle{empty}
\pagestyle{empty}

% using IEEEtran.bst settings
%\bstctlcite{IEEE_BSTcontrol}

\begin{abstract}
An explicit Model Predictive Control algorithm for large-scale structured linear systems is presented. We base our results on  Distributed and Localized Model Predictive Control (DLMPC), a  \emph{closed-loop} model predictive control scheme  based on the System Level Synthesis (SLS) framework wherein only local state and model information needs to be exchanged between subsystems for the computation and implementation of control actions. We provide an explicit solution for each of the subproblems resulting from the distributed MPC scheme. We show that given the separability of the problem, the explicit solution is only divided into three regions per state and input instantiation, making the point location problem very efficient. Moreover, given the locality constraints, the subproblems are of much smaller dimension than the full problem, which  significantly reduces the computational overhead of explicit solutions. We conclude with numerical simulations to demonstrate the computational advantages of our method, in which we show a large improvement in runtime per MPC iteration as compared with the results of computing the optimization with a solver online. 
\end{abstract}

\section{Introduction}\label{sec:introduction}

% History of the field
Model Predictive Control (MPC) has been shown to provide  solutions for many industrial applications, but its applicability was long limited to slow processes, since solving an optimal control problem online imposes a large computational burden. Explicit MPC was developed to overcome this issue, shifting the burden offline and reducing online computation to providing the evaluation of a piecewise function by relying on the principles of multiparametric programming \cite{bemporad_explicit_2000,bemporad_model_2000,pistikopoulos_-line_2002}. Despite its profound success, two important limitations restrict the applicability of MPC to large networks. On the one hand, explicit MPC has a limitation that is inherited from the computational complexity of multiparametric programming: finding a (piecewise) closed-form solution to an optimization problem becomes intractable for even modestly sized problems. On the other, even in the cases where the offline computation can be carried out, the solution is typically too complex to be evaluated efficiently online, in terms of both memory and runtime evaluation. 

% Related prior work
These problems in applying MPC to large networks relate to the fact that, in the worst case, complexity increases exponentially with the number of constraints \cite{morari_survey_2009}. Efforts have been made to circumvent these issues. First, the complexity of the offline computation  has been addressed by simplifying the MPC setup, using for example minimum-time formulations \cite{grieder_stabilizing_2005} or model reduction \cite{hovland_explicit_2008}, among others. Secondly, efforts have also been made to tackle the online limitations, i.e. to facilitate efficient solutions to the point-location problem; we elaborate further on this point in \sec{problem}. Examples of this are the partial enumeration method, where explicit MPC is implemented on only a subset of the constraints \cite{pannocchia_fast_2007}; and methods to optimally merge regions in order to reduce the number of partitions, as in\cite{geyer_optimal_2008}, \cite{holaza_nearly_2015} and the references therein. Once again, however, these methods are limited to systems of modest sizes and induce suboptimality in systems of large dimensions. Recent work has involved first formulating the network control problem as a distributed MPC problem, and then applying explicit MPC to the subproblems. Examples can be found in \cite{koehler_building_2013} and \cite{zheng_networked_2013}, both of which are developed for a specific application of systems with decoupled dynamics, where heuristics specific to the applications under consideration are developed. Although these approaches work well for the intended application, they do not generalize well.  In this setting, moreover, having an explicit solution can, perhaps counterintuitively, lead to prohibitive increases in runtime because distributed MPC approaches are usually iterative and require solving multiple optimization problems at each MPC iteration. Strategies more theoretical in focus  were developed in \cite{zarate_florez_explicit_2011}, and rely on a hierarchical structure with a global coordinator to tackle the distributed MPC problem. Purely distributed settings were analyzed in \cite{zhong_cooperation-based_2013}, and later in \cite{jiang_parallel_2019}. Although recursive feasibility and stability guarantees are provided, they rely on approximations, and as such, their solutions  are suboptimal.

% Limitations of prior work
Despite these efforts, the problem of how to make explicit MPC scalable, optimal, and applicable to large network settings remains an open question. This is especially relevant in distributed settings, where each subcontroller is typically repeatedly solving an optimization problem at every time step. Therefore, having an explicit solution in this regime will result in significant runtime improvements, since replacing each optimization problem with an explicit solution reduces total runtime by a factor of the number of iterations needed by each MPC subroutine.

% Purpose & summary of this work
In order to address this gap, we propose an explicit MPC solution that is applicable to large networks. We leverage the System Level Synthesis (SLS) \cite{wang_separable_2018,anderson_system_2019, chen_system_2019} framework, which provides a transparent and tractable way of dealing with distributed control synthesis problems. In particular, we rely on a novel parameterization of distributed closed-loop MPC policies introduced in \cite{amo_alonso_distributed_2019} such that the resulting synthesis problem is \emph{convex, structured, and separable}. Moreover, this parametrization allows for the MPC problem to be carried out over \emph{closed-loop policies} as opposed to open-loop inputs. We use the results of \cite{amo_alonso_distributed_2019}, which show that by exploiting the locality constraints allowed by the SLS parametrization, together with some mild separability assumptions, the MPC problem can be separated into small subproblems, each of which can solved in parallel by a different subsystem in the network. 
The main contributions of this work are to:
\begin{itemize}
\item Provide explicit solutions to the optimization problem each subcontroller in the system must solve in order to implement a distributed MPC algorithm.
\item Show that the explicit solution requires just 3 partitions of the solution space per system state/input instantiation, thus making the point-location problem trivial when solving for each of the instantiations sequentially.
\item Show, through simulations, that the complexity of the subproblems solved at each subsystem scales as $O(1)$ relative to the full size of the system while achieving optimal performance.
\end{itemize}

\subsubsection*{\textbf{Notation}} 
Bracketed indices  denote the time of the true system, i.e., the system input is  $u(t)$ at time $t$, not to be confused with $x_t$ which denotes the predicted state $x$ at time $t$. Superscripted variables, e.g. $x^k$,  correspond to the value of $x$ at the $k^{th}$ iteration of a given algorithm. To denote subsystem variables, we use square bracket notation, i.e., $[x]_{i}$ denotes the components of $x$ corresponding to subsystem $i$.  Calligraphic letters such as $\mathcal{S}$ (with the exception of $\mathcal A$ and $\mathcal B$) denote sets, and lowercase script letters such as $\mathfrak{c}$ denotes a subset of $\mathbb{Z}^{+}$, i.e., $\mathfrak{c}=\left\{1,...,n\right\}\subset\mathbb{Z}^{+}$.  Boldface lower and upper case letters such as $\mathbf{x}$ and $\mathbf{K}$ denote finite horizon signals and lower block triangular (causal) operators respectively (their form will be discussed in the relevant section).

  $\mathbf{K}(\mathfrak{r},\mathfrak{c})$ denotes the submatrix of $\mathbf{K}$ composed of the rows and columns specified by $\mathfrak{r}$ and $\mathfrak{c}$ respectively. The $\dagger$ superscript denotes the pseudo-inverse of a matrix.

\section{Problem Formulation} \label{sec:problem}

% DESCRIPTION OF MPC
We consider discrete-time systems of the form:
\begin{equation} \label{eq: LTI system}
x(t+1) = Ax(t)+Bu(t),\quad t = 0,1,\hdots, T,
\end{equation}

where $x(t)\in\mathbb{R}^{n}$ is the state, and $u(t)\in\mathbb{R}^{p}$ is the control input. 
The control input will be designed by an MPC scheme, where at each time $t$ the controller will solve the following optimal control problem: 
\begin{equation} \label{eq: MPC}
\begin{aligned}
& \underset{u_{t}, \gamma_t}{\text{min}} &  &\sum_{t=0}^{T-1}x_{t}^{\mathsf{T}}Q_{t}x_{t}+u_{t}^{\mathsf{T}}R_{t}u_{t}+x_{T}^{\mathsf{T}}Q_{T}x_{T}\\
& \ \text{s.t.} &  &\begin{aligned} 
    & x_{0} = x(t),\, \ x_{t+1} = Ax_{t}+Bu_{t},\, \ t=0,...,T-1,\\
    & x_{T}\in\mathcal{X}_{T},\, x_{t}\in\mathcal{X}_{t},\, u_{t}\in\mathcal{U}_{t} \, \ t=0,...,T-1, \\
    &u_{t} = \gamma_t(x_{0:t},u_{0:t-1}),
\end{aligned}
\end{aligned}
\end{equation}
where the cost matrices $Q_{t}$ and $R_{t}$ are positive semidefinite and positive definite respectively, $\mathcal{X}_t$ and $\mathcal{U}_t$ are polytopes containing the origin, and $\gamma_t(\cdot)$ are measurable functions of their arguments. 

The goal of this work is to design a strategy that allows for problem \eqref{eq: MPC} to have an \emph{explicit} solution, and that the computation time of such solution remains low and independent of the size of the system \eqref{eq: LTI system}. To do so we will exploit the underlying structure of the network and rely on a Distributed and Localized MPC (DLMPC) scheme \cite{amo_alonso_distributed_2019}, which ensures local computation and implementation of the closed-loop control policy. 

We view \eqref{eq: LTI system} as a network of $N$ interconnected subsystems. The state and control vectors  can be suitably partitioned as $[x]_i$ and $[u]_i$,  inducing a compatible block structure $[A]_{ij}$, $[B]_{ij}$ in the dynamics matrices $(A,B)$. In particular, this interconnection topology can be modeled as a time-invariant unweighted directed graph $\mathcal{G}_{(A,B)}(E,V)$, where each subsystem $i$ is identified with a vertex $v_{i}\in V$ and an edge $e_{ij}\in E$ exists whenever $[A]_{ij}\neq 0$ or $[B]_{ij}\neq 0$.

We further assume that the information exchange topology between subcontrollers  matches that of the underlying system, and can thus be modeled by the same graph $\mathcal{G}_{(A,B)}(E,V)$. Since we want the MPC control policy to respect the structure of the system, we impose that information exchange be localized to a subset of neighboring subcontrollers. To do so, we impose $d$-local information exchange constraints \cite{wang_localized_2016}, i.e., we require each subcontroller to only exchange their state and control input with neighboring controllers at most $d$-hops away according to the communication topology  $\mathcal{G}_{(A,B)}(E,V)$. This notion is formalized in the Definition \ref{def: in_out set} and an example is provided in Appendix B.
\begin{definition}{\label{def: in_out set}}
For a graph $\mathcal{G}(V,E)$, the \textit{d-outgoing set} of subsystem $i$ is 
\begin{equation*}
\textbf{out}_{i}(d) := \left\{v_{j}\ |\  \textbf{dist}(v_{i} \rightarrow v_{j} ) \leq d\in\mathbb{N} \right\}.
\end{equation*} The \textit{d-incoming set} of subsystem $i$ is 
\begin{equation*}
\textbf{in}_{i}(d) := \left\{v_{j}\ |\ \textbf{dist}(v_{j} \rightarrow v_{i} ) \leq d\in\mathbb{N} \right\}.
\end{equation*}
 Note that $v_i \in \textbf{out}_{i}(d)\cap \textbf{in}_{i}(d)$ for all $d\geq 0$. 
\end{definition}
In light of this interconnection topology, we want to construct a distributed MPC control algorithm such that both the synthesis and the implementation of the control input at each subsystem is \emph{localized}; that is, each closed-loop control policy at subsystem $i$ can be computed using only the state, control input and plant model from $d$-hop incoming neighbors $\textbf{in}_{i}(d)$, where $d$ - the size of the local neighborhood - is a design parameter. This can be achieved by imposing a $d$-local information exchange constraint on problem \eqref{eq: MPC}, so that each local control input takes the form: 
\begin{align}    \label{eq: info_constraints}
    [u_t]_i = \gamma_{i,t}\left([x_{0:t}]_{j\in \textbf{in}_i(d)},[u_{0:t-1}]_{j\in \textbf{in}_i(d)},
    [A]_{j,k \in \textbf{in}_i(d)},[B]_{j,k\in  \textbf{in}_i(d)}\right), 
\end{align}
 for all $t=0,\dots,T$ and $i=1,\dots,N$, where $\gamma_{i,t}$ is a measurable function of its arguments. 

As shown in  \cite{amo_alonso_distributed_2019}, under suitable structural compatibility assumptions between the cost function, state and input constraint sets, and information exchange constraints, the DLMPC algorithm introduced allows for both implementation and synthesis of a closed-loop distributed MPC control law in a localized manner by solving problem \eqref{eq: MPC} subject to the constraint \eqref{eq: info_constraints} for  the \emph{closed-loop system responses} of the system. 

% DESCRIPTION OF EXPLICIT MPC
The MPC algorithm proposed in \cite{amo_alonso_distributed_2019} is iterative, so an optimization problem is solved repeatedly until convergence. Simulations suggest that the cases where a closed-form solution exists enjoy $10\times$ faster convergence. 
Although constrained optimization problems rarely have a closed form, the explicit MPC approach allows for a piecewise solution to be computed offline.

Explicit MPC is based on the observation that at each time step, optimization \eqref{eq: MPC} remains constant except for the update in $x_{0}$. This $x_{0}$ can be seen as a parameter, and its value fully determines the solution to the optimization, i.e., $u^{*}(x_0)$. Explicit MPC makes the dependence between $u^{*}$ and $x_0$ explicit rather than it being implicitly obtained through the optimization problem. By leveraging the KKT conditions \cite{boyd_convex_2004}, one can solve  \eqref{eq: MPC} offline for all $x_0\in\mathcal{X}$ \cite{bemporad_explicit_2000} . The solution for each MPC subroutine (a convex quadratic program) is of the form
\begin{equation}
u^{*}(x_0) = \left\{\begin{array}{cc}
F_1x_0+g_1 & \text{if } H_1x_0 \leq h_1 \\
\vdots & \vdots \\ 
F_mx_0+g_m & \text{if } H_mx_0 \leq h_m,
\end{array}\right.
\label{eq: eMPC}
\end{equation} 
where $F_i,H_i$ and $f_i,h_i$ for $i=1,\dots,n$ are matrices and vectors respectively of the appropriate dimension (see \cite{bemporad_explicit_2000}). Given this explicit relation, the online optimization can be replaced by a \emph{point location problem}, i.e., the online problem boils down to i) finding which one of the constraints in expression \eqref{eq: eMPC} is satisfied given the current initial condition $x_0$, and ii) applying the appropriate control policy.

% WHY EXPLICIT MPC SUPERSEDES MPC, AND WHAT IS THE PROBLEM FOR LARGE NETWORKS
In order to solve an explicit MPC problem, two steps are required: (a) \emph{the offline step}: finding piecewise solution to the optimization problem \eqref{eq: MPC} given $x_0$, (b) \emph{the online step}: solving the point location problem. Clearly most of the computational effort takes place offline, reducing the online computational requirement. Explicit MPC has been shown to supersede conventional MPC in runtime while maintaining the same optimality and feasibility guarantees \cite{pistikopoulos_-line_2002}. However, this is only true for systems of small dimension, in which the number of constraints $m$ is a small. For systems with a large number of constraints, explicitly solving the offline optimization problem is prohibitive.

We address this problem and provide an explicit solution to the MPC problem \eqref{eq: MPC} subject to the constraint \eqref{eq: info_constraints} for arbitrarily large structured networks and using the DLMPC algorithm introduced in \cite{amo_alonso_distributed_2019}, which we summarize in the next section.

\section{Distributed and Localized MPC}\label{sec:dlmpc}
Consider the dynamical system \eqref{eq: LTI system} with an additional additive noise term $w(t)$, over a finite horizon $t = 0,...,T$. Applying a time varying state feedback control law $u_t= K_t(x_0,\hdots, x_t)$, where $K_t$ is a linear map to be designed. The closed-loop dynamics can be compactly written as
\begin{equation}\label{eq:SLScompact} 
\mathbf{x} = Z(\mathcal A+ \mathcal B\mathbf{K})\mathbf{x+\mathbf{w}},
\end{equation}
where $Z$ is the block-downshift matrix\footnote{A matrix with identity matrices along its first block sub-diagonal and zeros elsewhere},  $\mathcal A:=\mathrm{blkdiag}(A,A,...,A,0)$, and $\mathcal B:=\mathrm{blkdiag}(B,B,...,B,0)$. The vectors $\mathbf{x}$, $\mathbf{u}$, $\mathbf{w}$, are the finite horizon signals corresponding to state, control input, and disturbance respectively, and $\mathbf K$ represents the block matrix operator for the causal linear time-varying state-feedback controller. Rewriting~\eqref{eq:SLScompact} and $\mathbf{u=Kx}$ we obtain
\begin{equation}\label{eq:response}
\begin{split}
\mathbf{x} & = (I-Z(\mathcal A+\mathcal B\mathbf{K}))^{-1}\mathbf{w} =: \mathbf\Phi_x \mathbf w,\\
\mathbf{u} & = \mathbf{K}(I-Z(\mathcal A+\mathcal B\mathbf{K}))^{-1}\mathbf{w} =: \mathbf\Phi_u \mathbf w.
\end{split}
\end{equation}
The pair $\{ \mathbf \Phi_x, \mathbf \Phi_u \}$ is referred to as the \emph{system response} and one realization of the controller is given by $\mathbf K = \mathbf \Phi_u \mathbf \Phi_x^{-1}$\cite{anderson_system_2019}. In the SLS framework, control synthesis is (equivalently) reformulated as an optimization problem over system responses $\{ \mathbf \Phi_x, \mathbf \Phi_u \}$. The central result of SLS states that the resulting synthesis problem is convex in the system response matrices. We provide the theorem that leads to this conclusion in the appendix, for full details and to see examples of modeling many standard control problems in the SLS setting, see~\cite{anderson_system_2019}.

For finite horizon systems, the system response $\mathbf \Phi_x$ (the closed-loop map from disturbance to state) and  state $\mathbf x$, take the form
\begin{align*} 
  &\underbrace{{{\left[\begin{array}{cccc}\Phi_{x,0}[0] & & & \\ \Phi_{x,1}[1] & \Phi_{x,1}[0] & & \\ \vdots & \ddots & \ddots & \\ \Phi_{x,T}[T] & \dots & \Phi_{x,T}[1] & \Phi_{x,T}[0] \end{array}\right]}}} ~\text{and}~ \underbrace{\left[\begin{array}{c} x_{0}\\x_{1}\\\vdots\\x_{T}\end{array}\right]}, \\ & \hspace{3.15cm}\mathbf \Phi_x \hspace{4.55cm}\mathbf x
\end{align*}
respectively, where $\Phi_{x,i}[j]$ corresponds to the system response $\Phi_{x}$ synthesized at time $i$ and applied to disturbance $w_j$.\footnote{System~\eqref{eq: LTI system} does not contain a noise term. However in the SLS framework, we treat a non-zero initial condition as the first term in the disturbance sequence. The rest of the sequence is set to zero. } 
The system response from disturbance to control law, $\mathbf{\Phi}_u$, takes the same form, with the only difference being that it  runs to time $T-1$, and so does the controller $\mathbf{K}$. The control input and the disturbance signals follow analogously to the state. We denote the $k^{th}$ block column of $\mathbf{\Phi}_x$ as $\mathbf{\Phi}_x\{k\}$ i.e. $\mathbf{\Phi}_x\left\{0\right\}:=[\Phi_{x,0}[0]^{\mathsf{T}}\ \dots\ \Phi_{x,T}[T]^{\mathsf{T}}]^{\mathsf{T}}$. For compactness, sometimes we write $\mathbf{\Phi}:=[\mathbf{\Phi}_x^{\mathsf{T}}\; \mathbf{\Phi}_u^{\mathsf{T}}]^{\mathsf{T}}$.

As described in \cite{anderson_system_2019}, the controller can be implemented by   $\mathbf{u}=\mathbf{\Phi}_u\mathbf{\hat{w}}, \mathbf{\hat{x}}=(I - \mathbf{\Phi}_x)\mathbf{\hat{w}}, \mathbf{\hat{w}}=\mathbf{x}-\mathbf{\hat{x}}$, where $\mathbf{\hat{x}}$ can be interpreted as a nominal state trajectory, and $\mathbf{\hat{w}}=Z\mathbf{w}$ is a reconstruction of the disturbance with unit step delay. The advantage of this implementation is that structure imposed on the system response $\{\mathbf{\Phi}_x, \mathbf{\Phi}_u \}$ is mirrored in the controller structure.

One of the main advantages of taking this approach to reformulate the MPC problem is that optimization is done over \emph{closed-loop policies} instead of open-loop inputs. This allows for the synthesis of distributed and structured control policies in a convex manner. A similar approach was taken in \cite{goulart_optimization_2006}. The main difference is that the SLS framework also allows for a localized and distributed optimization of such policies as we will show.
\begin{lemma}
The MPC subproblem \eqref{eq: MPC} for system \eqref{eq: LTI system}  can be equivalently formulated in the SLS framework as:
\begin{equation}\label{eq: SLS MPC}
\begin{array}{rl}
\underset{\mathbf{\Phi}_x\left\{0\right\},\mathbf{\Phi}_u\left\{0\right\}}{\text{min}} & \left\Vert [C\ D] \begin{bmatrix}\mathbf{\Phi}_x\left\{0\right\}x_{0} \\ \mathbf{\Phi}_u\left\{0\right\}x_{0}\end{bmatrix}\right\Vert_F^2\\
 \text{s.t.} &  Z_{AB}\mathbf{\Phi}\left\{0\right\}=I,\, x_0 = x(t),\\
 	        & \mathbf{\Phi}_x\left\{0\right\}x_{0}\in \mathcal{X}^T, \mathbf{\Phi}_u\left\{0\right\}x_{0} \in\mathcal{U}^T,
\end{array}
\end{equation}
where we use $ Z_{AB}\mathbf{\Phi}\left\{0\right\}=I$ to compactly denote constraint \eqref{eq: constraint} in Theorem \ref{thm: SLS} (Appendix A), and $C$ and $D$ are constructed by arranging $Q_t^{\frac{1}{2}}$ and $R_t^{\frac{1}{2}}$ for all $t=1,\dots,T$ respectively in a block diagonal form. The constraints are encoded as   $\mathcal{X}^T := \big(\otimes_{t=0}^{T-1} \mathcal{X}\big) \otimes \mathcal{X}_T$, and similar for $\mathcal{U}^T$. 
\label{lem: MPC reformulation}
\end{lemma}
{\begin{proof} The reader is referred to $\S3$ of \cite{amo_alonso_distributed_2019}.\end{proof}}
As shown in~\cite{wang_separable_2018}, the cost function in problem \eqref{eq: MPC} encodes for the $H_2$-norm of the system responses.
For the remainder of the paper, we will overload notation and write $\mathbf{\Phi}_x$ and  $\mathbf{\Phi}_u$ in place of $\mathbf{\Phi}_x\left\{0\right\}$ and $\mathbf{\Phi}_u\left\{0\right\}$, given that no driving noise is present, only the first block columns of the system responses need to be computed.

The fact that problem \eqref{eq: SLS MPC} is solved over \emph{closed-loop policies} has many important implications, one of them being that any structure imposed on the system responses translates directly into the structure of the closed-loop map. This is true for local communication constraints \eqref{eq: info_constraints}, which can be transparently applied through \emph{locality constraints}.

\begin{definition}{\label{def: locality}}
Let $[\mathbf{\Phi}_{x}]_{ij}$ be the submatrix of system response $\mathbf{\Phi}_x$ describing the map from disturbance $[w]_{j}$ to the state $[x]_i$ of subsystem $i$. The map $\mathbf{\Phi}_{x}$ is \textit{d-localized} if and only if for every subsystem $j$, $[\mathbf{\Phi}_{x}]_{ij}=0\ \forall\ i\not\in\textbf{out}_{j}(d)$. The definition for \textit{d-localized} $\mathbf{\Phi}_u$ is analogous but with perturbations to control action $[u]_i$ at subsystem $i$.
\end{definition}

The constraint that  $\mathbf{\Phi}_{x}$ and $\mathbf{\Phi}_{u}$ are $d$-localized means that each subsystem only needs to collect information from its $d$-incoming set to implement the control law. Similarly, it only needs to share information with its $d$-outgoing set to allow for other subsystems to implement their respective control laws. Notice that $d$-localized system responses are system responses with suitable sparsity patterns. 

\begin{definition}{\label{def: locality constraints}}
A subspace $\mathcal{L}_d$ enforces a $d$\textit{-locality constraint} if $\mathbf{\Phi}_{x},\mathbf{\Phi}_{u}\in\mathcal{L}_d$ implies that $\mathbf{\Phi}_{x}$ is $d$-localized and $\mathbf{\Phi}_{u}$ is $(d+1)$-localized.  \end{definition}

When locality constraints are introduced together with the  following compatibility assumptions between the cost function and state and input constraints, $d$-local information exchange constraints allow for a distributed and localized synthesis of the MPC problem.
\begin{assumption}{\label{assump: locality}}
Matrices $Q_{t}$ and $R_t$ for all $t=1,\dots,T$ in formulation (\ref{eq: MPC}) are structured such that $x_{t}^{\mathsf{T}}Q_{t}x_{t} = \sum_{i=1}^N [x_{t}]_{i}^{\mathsf{T}}[Q_{t}]_i[x_{t}]_{i}$, and $u_{t}^{\mathsf{T}}R_{t}u_{t} =\sum_{i=1}^N [u_{t}]_{i}^{\mathsf{T}}[R_{t}]_i[u_{t}]_{i}$. The constraint sets in formulation (\ref{eq: MPC}) are such that $x\in\mathcal{X}=\mathcal{X}_1\times ... \times \mathcal{X}_n$, where $x \in \mathcal{X}$ if and only if $[x]_{i}\in\mathcal{X}_{i}$ for all $i$, and idem for $\mathcal{U}$. 
\end{assumption}

We can now formulate the DLMPC subproblem by suitably incorporating locality constraints as well as Assumption \ref{assump: locality} into the SLS based MPC subproblem \eqref{eq: SLS MPC}.
\begin{equation}\label{eq: SLS MPC local}
\begin{array}{rl}
\underset{\mathbf{\Phi}_{x},\mathbf{\Phi}_{u}}{\text{min}} & \sum_{i=1}^N \left\Vert [[C]_i\ [D]_i] \begin{bmatrix}[\mathbf{\Phi}_xx_{0} ]_i\\ [\mathbf{\Phi}_ux_{0}]_i\end{bmatrix}\right\Vert_F^2\\
\text{s.t.} & Z_{AB}\mathbf{\Phi}=I,\, x_0 = x(t),\\ & [\mathbf{\Phi}_{x}x_0]_{i}\in\mathcal{X}_{i},\ [\mathbf{\Phi}_{u}x_0]_{i}\in\mathcal{U}_{i}, \\& i=1,\dots,N,\, \mathbf{\Phi}_{x},\mathbf{\Phi}_{u}\in\mathcal{L}_{d},
\end{array}
\end{equation}
where $[C]_i$ and $[D]_i$ are defined so as to be compatible with the $[Q_{t}]_i$ and $[R_{t}]_i$ defined in Assumption \ref{assump: locality}.

One can now exploit the separability\footnote{For details on row-wise and column-wise separability the reader is referred to \cite{amo_alonso_distributed_2019} and references therein.} of the problem using an algorithm for distributed optimization such as the Alternating Direction Method of the Multipliers (ADMM) \cite{boyd_distributed_2010}. In this way, each subsystem can solve for a reduced problem obtained as a subproblem of problem \eqref{eq: SLS MPC local} (see \cite{amo_alonso_distributed_2019}). Moreover, since the system responses are restricted to be $d$-localized, i.e., that $\mathbf{\Phi}_{x},\mathbf{\Phi}_{u}\in\mathcal{L}_d$, the resulting subproblem variables are sparse by Definition \ref{def: locality}, which allows for a significant reduction in the dimension of each local subproblem solved by each subsystem. 
In order to highlight the decomposable nature of the solution,  we require some additional notation.  We say that $\mathbf{\Phi}x_0\in\mathcal{P}$ iff  $\mathbf{\Phi}_{x}x_0 \in\mathcal{X}^T$\ \emph{and} $\mathbf{\Phi}_{u}x_0\in\mathcal{U}^T$, and let $[\mathcal P]_i$ denote the appropriate local subset. Similarly, we define $[\hat{C}]_i:=[[C]_i\ [D]_i]$. Define $[\mathbf{\Phi}]_{i_{r}}:=\mathbf{\Phi}(\mathfrak{s_{\mathfrak{r}_{i}}},\mathfrak{r}_{i})$ and $[\mathbf{\Phi}]_{i_{c}}:=\mathbf{\Phi}(\mathfrak{c}_{i},\mathfrak{s_{\mathfrak{c}_{i}}})$, where the sets  $\mathfrak{r}_{i}$ and  $\mathfrak{c}_{i}$ correspond to the rows and columns that controller $i$ is solving for, and the set $\mathfrak{s_{\mathfrak{r}_{i}}}~ (\mathfrak{s}_{\mathfrak{c}_{i}})$ is the set of columns (rows) associated to the rows (columns) in $\mathfrak{r}_{i}$ ($\mathfrak{c}_{i}$) by the locality constraints $\mathcal{L}_{d}$.

The ADMM problem solved (in parallel) by each subcontroller $i$ is:
 \begin{subequations}\label{eq: MPC ADMM localized 1}
\begin{align}
& [\mathbf{\Phi}]_{i_{r}}^{k+1} = \nonumber\\
& \hspace{-0.2cm} \left\{
\begin{aligned}
&\underset{\mathbf{[\mathbf{\Phi}]}_{i_{r}}}{\text{argmin}} && \hspace{-0.4cm} \left\Vert [\hat{C}]_i [\mathbf{\Phi}]_{i_r}[x_{0} ]_{i_r}\right\Vert_F^2+\frac{\rho}{2}\left\Vert[\mathbf{\Phi}]_{i_{r}}-[\mathbf{\Psi}]_{i_{r}}^{k}+[\mathbf{\Lambda}]_{i_{r}}^{k}\right\Vert^{2}_{F}\\
&\text{s.t.} && \begin{aligned}
     &[\mathbf{\Phi}]_{i_{r}}[x_0]_{i_{r}}\in[\mathcal P]_i\\
\end{aligned}
\end{aligned}\right\} \label{eq: MPC ADMM localized 1 - row}
\\
& \begin{aligned}
& [\mathbf{\Psi}]_{i_{c}}^{k+1} & = \big([\mathbf{\Phi}]_{i_{c}}^{k+1}+[\mathbf{\Lambda}]_{i_{c}}^{k}\big)+[Z_{AB}]_{i_{c}}^{\dagger}\Big([I]_{i_{c}}-
- [Z_{AB}]_{i_{c}}\big([\mathbf{\Phi}]_{i_{c}}^{k+1}+[\mathbf{\Lambda}]_{i_{c}}^{k}\big)\Big), 
\end{aligned} \label{eq: MPC ADMM localized 1 - column}
\\
& [\mathbf{\Lambda}]_{i_{r}}^{k+1} =[\mathbf{\Lambda}]_{i_{r}}^{k}+[\mathbf{\Phi}]_{i_{r}}^{k+1}-[\mathbf{\Psi}]_{i_{r}}^{k+1}. \label{eq: MPC ADMM localized 1 - lagrange} 
\end{align}
\end{subequations}

When considering the row-wise subproblem \eqref{eq: MPC ADMM localized 1 - row} evaluated at subsystem $i$, the $j^{\text{th}}$ column of the $i^{\text{th}}$ subsystem row partition $\mathbf{\Phi}_x(\mathfrak{r}_{i},:)$ and $\mathbf{\Phi}_u(\mathfrak{r}_{i},:)$ is nonzero only if  $j\in\textbf{in}_j(d)$ and $j\in\textbf{in}_j(d+1)$, respectively.  It follows that subsystem $i$ only requires a corresponding subset of the local sub-matrices $[A]_{k,\ell}, [B]_{k,\ell}$ to solve its respective subproblem. All column/row/matrix subsets described above can be found algorithmically as described in Appendix A of \cite{wang_separable_2018}.  

Algorithm \ref{alg: I} summarizes the implementation at subsystem $i$ of the ADMM-based solution to the DLMPC subproblem (\ref{eq: SLS MPC local}) under assumption \ref{assump: locality}. This algorithm is run in parallel by each subcontroller. 

\setlength{\textfloatsep}{0pt}% Remove \textfloatsep
\begin{algorithm}[h]
\caption{Subsystem $i$ DLMPC implementation}\label{alg: I}
\begin{algorithmic}[1]
\State \textbf{input:} convergence tolerance parameters $\epsilon_p>0$, $\epsilon_d>0$
\State Measure local state $[x(t)]_{i}$.
\State Share the measurement with $\textbf{out}_{i}(d)$.
\State Solve optimization problem (\ref{eq: MPC ADMM localized 1 - row}).
\State Share $[\mathbf{\Phi}]_{i_{r}}^{k+1}$ with $\textbf{out}_{i}(d)$. Receive the corresponding $[\mathbf{\Phi}]_{j_{r}}^{k+1}$ from $\textbf{in}_{i}(d)$ and build $[\mathbf{\Phi}]_{i_{c}}^{k+1}$.
\State Solve optimization problem (\ref{eq: MPC ADMM localized 1 - column}) via the closed form solution (\ref{eq: MPC ADMM localized 1 - column}).
\State Share $[\mathbf{\Psi}]_{i_{c}}^{k+1}$ with $\textbf{out}_{i}(d)$. Receive the corresponding $[\mathbf{\Phi}]_{j_{c}}^{k+1}$ from $\textbf{in}_{i}(d)$ and build $[\mathbf{\Psi}]_{i_{r}}^{k+1}$.
\State Perform the multiplier update step (\ref{eq: MPC ADMM localized 1 - lagrange}).
\State Check convergence as $\left\Vert[\mathbf{\Phi}]_{i_{r}}^{k+1}-[\mathbf{\Psi}]_{i_{r}}^{k+1}\right\Vert_F\leq\epsilon_{p}$ and $\left\Vert[\mathbf{\Psi}]_{i_{r}}^{k+1}-[\mathbf{\Psi}]_{i_{r}}^{k}\right\Vert_F\leq\epsilon_{d}$. 
\State If converged, apply computed control action $[u_0]_i = [\Phi_{u,0}[0]]_{i_{r}}[x_0]_{\mathfrak s_{\mathfrak{r_i}}}$\footnote{$[x_0]_{\mathfrak s_{\mathfrak{r_i}}}$ denotes the subset of elements of $x_0$ associated with the columns in $\mathfrak s_{\mathfrak{r_i}}$, such that $[\Phi_u^{0,0}]_{i_{r}}x_0 =[\Phi_u^{0,0}]_{i_{r}}[x_0]_{\mathfrak s_{\mathfrak{r_i}}}$}, and return to 2, otherwise return to 4. 
\end{algorithmic}
\end{algorithm}

Notice that despite the fact that a DLMPC controller can be synthesized locally  - so the subproblems are of small dimension - step 4 requires solving an optimization problem online, which significantly slows down runtime. In the next section, we will illustrate how we can use Algorithm \ref{alg: I} as a baseline for a distributed and localized \emph{explicit solution} to MPC that scales gracefully with the size of the system.

\section{Explicit MPC}\label{sec:explicit_mpc}
In the previous section we saw that Algorithm 1~\cite{amo_alonso_distributed_2019} allows us to solve the MPC problem \eqref{eq: MPC} in a distributed and localized manner:  each of the subsystems solves a subproblem of small dimension, and coordinates locally with its neighbors by means of an iterative scheme. This algorithm requires each subsystem to solve an optimization problem online, several times per MPC iteration (since each MPC iteration requires several iterations among subsystems as illustrated in Algorithm \ref{alg: I}). This motivates the need for an explicit solution that will reduce this computation burden. In particular, since subproblems \eqref{eq: MPC ADMM localized 1 - column} and \eqref{eq: MPC ADMM localized 1 - lagrange} can be solved in closed form, it is only subproblem \eqref{eq: MPC ADMM localized 1 - row} that will require an explicit solution. To achieve this, we take inspiration from the procedure used by the influential paper  \cite{bemporad_explicit_2000} and perform a similar analysis. However, the structure of our problem differs significantly from the original formulation in \cite{bemporad_explicit_2000}, both in the optimization variables and in the way the parameters enter the problem. Hence, the original derivation is not directly applicable to this case. 
We proceed as follows: assume $\rho>0$ is a  constant scalar. Let us define
\begin{equation*}
M:=\big(2x_0x_0^\intercal+\rho I\big)^{-1}.
\end{equation*}
Note that the structure of the matrix $2x_0x_0^\intercal+\rho I $ allows us to compute its inverse, $M$, very cheaply using the Sherman--Morrison formula.
\begin{lemma}\label{lem: explicit}
Let  $\Phi$ and $a$ be row vectors, $x_0$ column vector of compatible dimension, and $b_1,b_2$ scalars. Then, the optimal solution to
\begin{equation}
\begin{aligned}
&\underset{\Phi}{\text{min}} &&  \left\vert\Phi x_0\right\vert+\frac{\rho}{2}\left\Vert\Phi-a\right\Vert^{2}_{2}\\
&\text{s.t.} && \begin{aligned}
     &b_2\leq\Phi x_0\leq b_1,
\end{aligned}
\end{aligned}
\label{eq: explicit optimization}
\end{equation}
is
\begin{equation}\label{eq: explicit}
 \Phi^{\star}=\big(\rho a - \lambda x_0^\intercal\big)\big(2x_0x_0^\intercal+\rho I\big)^{-1},
\end{equation}
where 
\begin{equation*}
\lambda = \left\{\begin{array}{cc}
\frac{\rho a M x_0-b_1}{x_0^\intercal M x_0} & \text{if } \rho a M x_0-b_1>0 \\
\frac{\rho a M x_0-b_2}{x_0^\intercal M x_0} & \text{if } \rho a M x_0-b_2<0 \\
0 & \text{otherwise}.
\end{array}\right.
\label{eq: lambda}
\end{equation*}
\end{lemma}
{\begin{proof}
Apply the KKT conditions to optimization \eqref{eq: explicit optimization}. In particular, the stationarity condition  is:
\begin{align*}
    \nabla_{\Phi}\Big(\left\vert\Phi^{\star}x_0\right\vert+ & \frac{\rho}{2}\left\Vert\Phi^{\star}-a\right\Vert^{2}_{2}\Big)+\lambda_1\nabla_{\Phi}\Big(\Phi^{\star}x_0 - b_1\Big)\\ &+\lambda_2\nabla_{\Phi}\Big(-\Phi^{\star}x_0 + b_2\Big)=0,
 \end{align*}
where $\lambda_1$ and $\lambda_2$ represent two scalar Lagrange multipliers whose values are unknown. This leads to the following result for the optimal $\Phi$ as a function of the unknown $\lambda_1$ and $\lambda_2$:
\begin{equation}
    \Phi^{\star}=\big(\rho a - (\lambda_1-\lambda_2)x_0^\intercal\big)\big(2x_0x_0^\intercal+\rho I\big)^{-1}.
    \label{eq: Phi optimal}
\end{equation}
Notice that by  Slater's condition (Chapter 5 in \cite{boyd_convex_2004}) strong duality holds for problem \eqref{eq: explicit optimization}. Hence, we can make use of the dual problem to find the optimal solution. The dual problem can be written as:
\begin{align*}
  \underset{\lambda_1,\lambda_2\geq0}{\max}\ \left\vert\Phi^{\star}x_0\right\vert+& \frac{\rho}{2}\left\Vert\Phi^{\star}-a\right\Vert^{2}_{2}-\lambda_1(b_1-\Phi^{\star} x_0)-\lambda_2(-b_2+\Phi^{\star} x_0).
  \end{align*}
After substituting   $\Phi^{\star}$ into the dual problem above, the cost function becomes a quadratic function of $\lambda:=[\lambda_1\ \lambda_2]^\intercal$. In particular, after some algebraic manipulations one can conclude that the dual problem is a quadratic program equivalent to:
\begin{equation}
  \underset{\lambda\geq0}{\max}\ \ \lambda^\intercal c_2 \lambda + c_1 \lambda,
  \label{eq: dual}
\end{equation}
where $c_2 = \frac{1}{2}x_0^\intercal M x_0 \begin{bmatrix}-1 & &1 \\ 1 & &-1\end{bmatrix}$ and $c_1 = [\rho a M x_0 - b_1\ \ -\rho a M x_0 + b_2] $. 

In order to compute the value of $\lambda$, we exploit complementary slackness:
\begin{equation*}
\lambda_1(\Phi x_0 - b_1)=0, \quad \text{and} \quad
\lambda_2(-\Phi x_0 + b_2)=0.
\end{equation*}
This condition makes evident that $\lambda_1$ and $\lambda_2$ cannot be both nonzero, since by assumption $b_1<b_2$. Hence, let us assume without loss of generality that $\lambda_2=0$. The solution to problem \eqref{eq: dual} for $\lambda_1$ is as follows:
\begin{equation*}
    \lambda_1=\left\{\begin{array}{cc} \frac{\rho a M x_0-b_1}{x_0^\intercal M x_0} &\text{if }\lambda_1>0,\\ 0 &\text{otherwise.} \end{array}\right.
\end{equation*}

The form for $\lambda_1=0$ follows a similar structure. Notice that the matrix $M$ is by definition positive definite. Hence, $x_0^\intercal M x_0>0$ for all $x_0\neq 0$ and the sign of $\lambda_1$ is purely determined by the sign of $a M x_0-b_1$. This allows us to define the closed form solution for $\lambda$, and therefore for $\Phi$,  in a piecewise manner depending on the region. The criteria are specified in Table \ref{tab: regions}.

\begin{table}[h]
\begin{center}
\begin{tabular}{|l|l|}
\hline
 \begin{tabular}{l}
    Region in which $x_0$ lies
  \end{tabular} &
  \begin{tabular}{l}
     Corresponding solution for $\lambda$
  \end{tabular} 
  \\ \hline \hline
    $\rho a M x_0-b_1>0$ &
    $\lambda_1=\frac{\rho a M x_0-b_1}{x_0^\intercal M x_0}$, $\lambda_2=0$
  \\ \hline
     $-\rho a M x_0+b_2>0$ & 
     $\lambda_1=0$, $\lambda_2=\frac{-\rho a M x_0+b_2}{x_0^\intercal M x_0}$
  \\ \hline
  \begin{tabular}{l}
    $\rho a M x_0-b_1<0,$\\ $-\rho a M x_0+b_2<0$
  \end{tabular}  &
    $\lambda_1=0$, $\lambda_2=0$
  \\ \hline
\end{tabular}
\caption{Partition of the space of $x_0$ into the different regions that lead to different solutions for $\lambda$.}
\label{tab: regions}
\end{center}
\end{table}

Recall that from optimization \eqref{eq: explicit optimization}, the problem is only feasible if $b_1<b_2$, hence the regions defined in Table \ref{tab: regions} are disjoint and well-defined. Leveraging the entries of Table \ref{tab: regions} and equation \eqref{eq: Phi optimal}, one can find the explicit solution \eqref{eq: explicit}. 
\end{proof}}

We  now apply  Lemma \ref{lem: explicit} to step 4 of Algorithm \ref{alg: I} i.e., solving problem \eqref{eq: MPC ADMM localized 1 - row}.  Notice that when only safety and saturation constraints are allowed, problem \eqref{eq: MPC ADMM localized 1 - row} is:
\begin{equation}
\begin{aligned}
&\underset{\mathbf{[\mathbf{\Phi}]}_{i_{r}}}{\text{min}} &&  \left\Vert [\hat{C}]_i [\mathbf{\Phi}]_{i_r}[x_0]_{\mathfrak s_{\mathfrak{r_i}}}\right\Vert_F^2+\frac{\rho}{2}\left\Vert[\mathbf{\Phi}]_{i_{r}}-[\mathbf{\Psi}]_{i_{r}}^{k}+[\mathbf{\Lambda}]_{i_{r}}^{k}\right\Vert^{2}_{F}\\
&\text{s.t.} && \begin{aligned}
     &\begin{bmatrix} [\mathbf x^{min}]_i\\ [\mathbf u^{min}]_i\end{bmatrix}\leq[\mathbf{\Phi}]_{i_{r}}[x_0]_{\mathfrak s_{\mathfrak{r_i}}}\leq\begin{bmatrix} [\mathbf x^{max}]_i\\ [\mathbf u^{max}]_i\end{bmatrix}.\\
\end{aligned} 
\end{aligned}
\label{eq: theorem}
\end{equation}

Given the separability properties of the Frobenius norm and the constraints, this optimization problem can further be separated into \emph{single rows} of $[\mathbf{\Phi}]_{i_r}$ and $[\mathbf{\Psi}]_{i_{r}}^{k}-[\mathbf{\Lambda}]_{i_{r}}^{k}$. Notice that this is true for the first term of the objective function as well, since $[\hat{C}]_i$ is a diagonal matrix  by Assumption \ref{assump: locality}, so its components can be treated as factors multiplying each of the rows accordingly.

It is important to note that by definition of $[\mathbf{\Phi}]_{i_r}$, 
$[\mathbf{x}]_i = [\mathbf{\Phi}]_{i_r}[x_0]_{\mathfrak{s}_{\mathfrak{r_i}}}$.
Hence, each row of $[\mathbf{\Phi}]_{i_r}$ multiplied with $[x_0]_{\mathfrak s_{\mathfrak{r_i}}}$ precisely corresponds to a given component of the a state/input instantiation, i.e., $[x_t]_i$ /$[u_t]_i$.

We can now consider one of the single-row subproblems resulting from this separation and rename its variables, where $\Phi$ represents the given row of $[\mathbf{\Phi}]_{i_{r}}$, $a$ represents the corresponding row of $[\mathbf{\Psi}]_{i_{r}}^{k}-[\mathbf{\Lambda}]_{i_{r}}^{k}$ and $x_0,\ b_1$ and $b_2$ correspond to the  elements of $[x_0]_{\mathfrak s_{\mathfrak{r_i}}}$, $[x^{min}_t]_i$/$[u^{min}_t]_{i}$ and $[x^{max}_t]_{i}$/$[u^{max}_t]_{i}$ respectively. Without loss of generality we set each nonzero component of $[\hat{C}]_i$ to be equal to 1. By noting that for the inner product $\Phi x_0$, it holds that $\left\Vert\Phi x_0\right\Vert_F^2=\left\vert\Phi x_0\right\vert$, and for any vector, the Frobenius norm is equivalent to the 2-norm, i.e., $\left\Vert\Phi-a\right\Vert^{2}_{F}=\left\Vert\Phi-a\right\Vert^{2}_{2}$, we can directly apply Lemma \ref{lem: explicit} to each of the single-row subproblems in which the problem \eqref{eq: theorem} can be separated.

Hence, by Lemma \ref{lem: explicit}, an explicit solution exists for optimization \eqref{eq: MPC ADMM localized 1 - row}. Thus, step 4 in Algorithm \ref{alg: I} can be solved explicitly. Notice that all other computation steps in Algorithm \ref{alg: I} have closed-form solutions. The following theorem follows naturally from the previous discussion.

\begin{theorem}
Given the MPC problem \eqref{eq: MPC} subject to the information constraint \eqref{eq: info_constraints}, and assume that the constraint sets $\mathcal{X}_t$ and $\mathcal{U}_t$ are of the form $x^{min}_t\leq x_t\leq x^{max}_t$ and $u^{min}_t\leq u_t\leq u^{max}_t$ for all $t=1,\dots,T$,  then~\eqref{eq: MPC} can be solved via the iterative distributed and localized explicit solution presented in Algorithm \ref{alg: I explicit}. 
\label{thm: explicit}
\end{theorem}

\begin{algorithm}[h]
\caption{Subsystem $i$ explicit DLMPC implementation}\label{alg: I explicit}
\begin{algorithmic}[1]
\State \textbf{input:} convergence tolerance parameters $\epsilon_p>0$, $\epsilon_d>0$
\State Measure local state $[x(t)]_{i}$.
\State Share the measurement with $\textbf{out}_{i}(d)$.
\State Solve for each of the rows of  $[\mathbf{\Phi}]_{i_{r}}^{k+1}$ sequentially via the explicit solution \eqref{eq: explicit} with the appropriate variable renaming.
\State Share $[\mathbf{\Phi}]_{i_{r}}^{k+1}$ with $\textbf{out}_{i}(d)$. Receive the corresponding $[\mathbf{\Phi}]_{j_{r}}^{k+1}$ from $\textbf{in}_{i}(d)$ and build $[\mathbf{\Phi}]_{i_{c}}^{k+1}$.
\State Compute $[\mathbf{\Psi}]_{i_{c}}^{k+1}$ via the closed form solution (\ref{eq: MPC ADMM localized 1 - column}).
\State Share $[\mathbf{\Psi}]_{i_{c}}^{k+1}$ with $\textbf{out}_{i}(d)$. Receive the corresponding $[\mathbf{\Phi}]_{j_{c}}^{k+1}$ from $\textbf{in}_{i}(d)$ and build $[\mathbf{\Psi}]_{i_{r}}^{k+1}$.
\State Perform the multiplier update step (\ref{eq: MPC ADMM localized 1 - lagrange}).
\State Check convergence as $\left\Vert[\mathbf{\Phi}]_{i_{c}}^{k+1}-[\mathbf{\Psi}]_{i_{c}}^{k+1}\right\Vert_F\leq\epsilon_{p}$ and $\left\Vert[\mathbf{\Psi}]_{i_{c}}^{k+1}-[\mathbf{\Psi}]_{i_{c}}^{k}\right\Vert_F\leq\epsilon_{d}$. 
\State If converged, apply computed control action $[u_0]_i = [\Phi_{u,0}[0]]_{i_{r}}[x_0]_{\mathfrak s_{\mathfrak{r_i}}}$,
and return to step 2, otherwise return to step 4. 
\end{algorithmic}
\end{algorithm}

Convergence can be shown through the standard  ADMM convergence results~\cite{boyd_distributed_2010}. Recursive feasibility and stability results can be found in \cite{amo_alonso_distributed_2019}.

\begin{remark}
Contrary to conventional explicit MPC (where regions are computed offline and the online problem reduces to a point location problem), our approach is to solve  step 4 in Algorithm \ref{alg: I explicit} explicitly. This ensures that all the steps in the algorithm  are solved in closed form or via an explicit solution, hence we refer to our algorithm as explicit MPC.
\end{remark}

\begin{figure}[t]
    \centering
    \includegraphics[width=0.99\columnwidth]{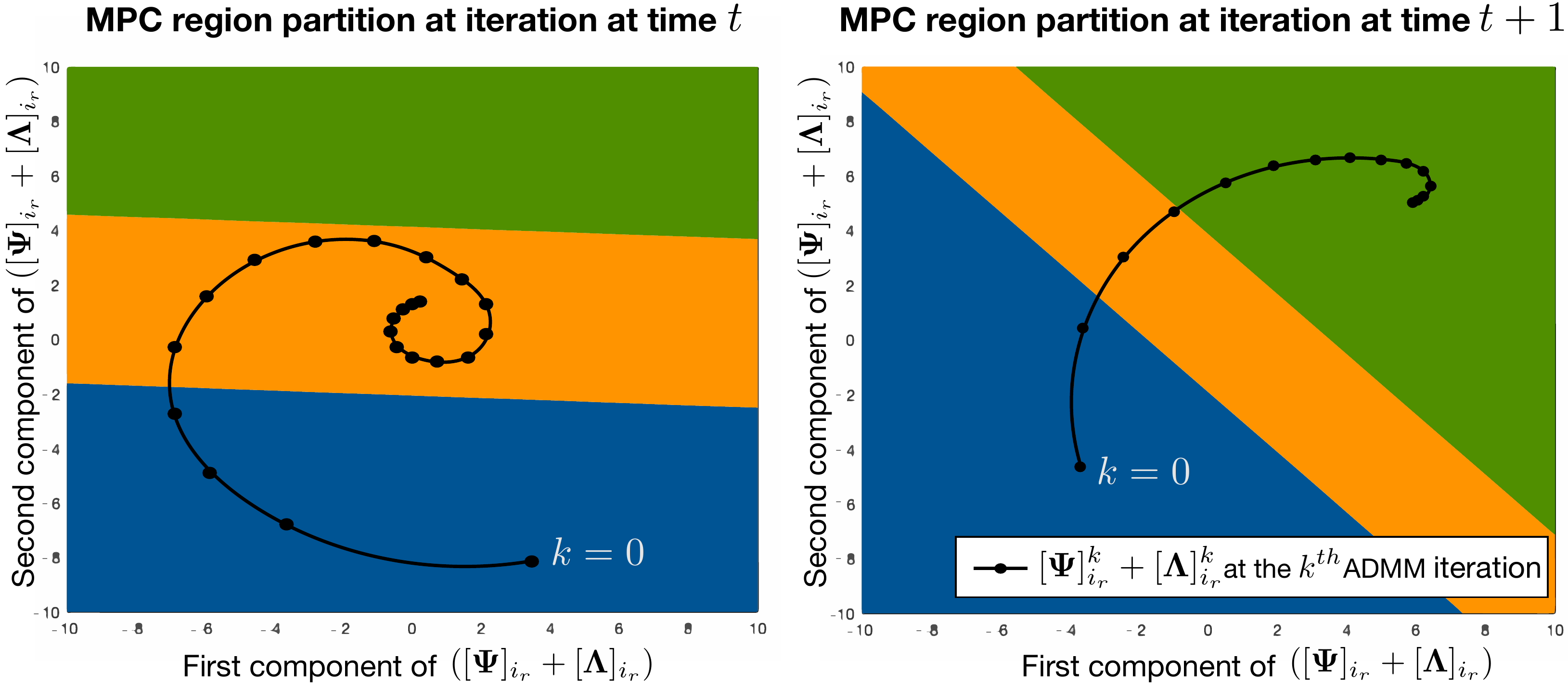}
    \caption{Illustration of the regions and parameter location over two MPC iterations, and the necessary ADMM iterations until convergence in each of the MPC iterations. For simplicity in the representation, we consider the parameters in twodimensions. \vspace{0.5cm}}
    \label{fig: regions}
\end{figure}

Another difference between standard explicit MPC and our formulation is that, in our case, the regions are not defined by polytopes of $x_0$ (as in~\eqref{eq: eMPC}). In our case, $M$ depends on $x_0$, thus, at each MPC iteration, new regions are computed as an explicit function of $x_0$, and for the subsequent ADMM iterations (within each MPC iteration) the parameter of the optimization problem is the corresponding row of $[\mathbf{\Psi}]_{i_{r}}^{k}-[\mathbf{\Lambda}]_{i_{r}}^{k}$ denoted as $a$ in \eqref{eq: explicit}. The regions are indeed affine with respect to this parameter.  Note that $x_0$ remains fixed within each MPC iteration. This idea is illustrated in Figure \ref{fig: regions}, where we illustrate the different regions involved in the computation of a given row of $[\mathbf{\Phi}]_{i_{r}}$ for two MPC iterations. In order to not overload notation, in this example we denote a single row of matrices $[\mathbf{\Phi}]_{i_{r}}$ and $[\mathbf{\Psi}]_{i_{r}}^{k}-[\mathbf{\Lambda}]_{i_{r}}^{k}$ with the same notation as the whole matrices themselves, i.e., $[\mathbf{\Phi}]_{i_{r}}$ and $[\mathbf{\Psi}]_{i_{r}}^{k}-[\mathbf{\Lambda}]_{i_{r}}^{k}$ denote a row of the homonymous matrices.

We now analyze the complexity of Algorithm \ref{alg: I explicit}. In particular, steps 3, 5 and 7 are only concerned with communication between the $d$-local neighbors. The computational complexity is determined by the update steps 4, 6 and 8. Similarly to Algorithm \ref{alg: I}, steps 6 and 8 boil down to the evaluation of a closed form expression involving matrix operations, so the complexity of these steps is $O(d^{2}T)$. Step 4 consists of a point location problem followed by a matrix multiplication. The complexity of solving for each row of $[\mathbf{\Phi}]_{i_{r}}$ results in $O(d^2)$, since the point location problem involves only 3 regions and the size of the matrices is $O(d^2)$. Given that each subsystem performs this operation sequentially for each of the rows in $[\mathbf{\Phi}]_{i_{r}}$, the complexity of step 4 is also  $O(d^2T)$. This is in contrast with Algorithm \ref{alg: I}, where step 4 consisted of solving an optimal problem with $O(d^{2}T)$ optimization variables and $O(dT)$ constraints. The significant overhead reduction given by the explicit solution in Algorithm \ref{alg: I explicit} is due to the simplicity of the point location problem. As stated in Lemma \ref{lem: explicit}, the space of the solution is partitioned into 3 regions per state/input instantiation, and this is independent of the size of the global system $N$, the size of the locality region $d$ and the total number of constraints. Hence, the complexity is dominated by the matrix multiplication needed to compute the explicit solution $\Phi^{\star}$. Notice that the locality constraints provide a computational advantage when $d\ll N$.

\section{Simulation experiments}\label{sec:simulations}

Consider a network with a chain topology where each node in the chain is a two-state system that evolves with dynamics
\begin{equation*}
[x(t+1)]_{i}=[A]_{ii}[x(t)]_{i}+\sum_{j\in\textbf{in}_{i}(d)}[A]_{ij}[x(t)]_{j}+[B]_{ii}[u(t)]_{i},
\end{equation*}
\vspace{-2mm}
where 
\begin{equation*}
[A]_{ii}=\begin{bmatrix}
   1  & 0.1 \\
  -0.3 &  0.7
\end{bmatrix}, \ [A]_{ij}=\begin{bmatrix}
   0  & 0 \\
   0.1 &  0.1
\end{bmatrix}, \ [B]_{ii}=\begin{bmatrix}
   0 \\
   0.1 
\end{bmatrix}.
\end{equation*}

The MPC  cost function is $$f(x,u) = \sum_{i=1}^{N} \sum_{t=1}^{T-1} \|[x(t)]_i\|_2^2+\|[u(t)]_i\|_2^2+\|[x(T)]_i\|_2^2,$$ where the time horizon is $T=5$, the locality parameter is set to $d=1$, and we vary the number of subsystems, $N$, considered. Specifically  $N\in\{10,50,100,200\}$.

We consider three control scenarios:
\begin{itemize}
\item \textbf{Case 1} (from \cite{amo_alonso_distributed_2019}): The \emph{unconstrained} MPC problem described above. It can be solved directly in closed form with Algorithm \ref{alg: I}. This will provide us with a baseline to establish the efficiency of our method.
\item \textbf{Case 2} (from \cite{amo_alonso_distributed_2019}): The MPC problem described above subject to constraints, solved using  Algorithm \ref{alg: I}. The (asymmetric) constraints considered are:
\begin{equation*}
-0.2\le[x(t)]_{i,1} \le1.2 \quad  \text{for } t=1,...,T,
\end{equation*}
where $[x]_{i,1}$ denotes the first state in the two-state subsystem $i$.
\item \textbf{Case 3}: Same setup at case 2, solved using   our explicit MPC scheme -- Algorithm \ref{alg: I explicit}.
\end{itemize}

This comparison illustrates the runtime improvement when using an explicit solution. We observe that the computational runtime of the explicit MPC scheme is much faster than its alternative solution with a solver, and it is close to being as fast as a closed-form solution due to the fact that there are only 3 regions per state/input instantiation. We also illustrate the relative error of the cost function when computed with Algorithm \ref{alg: I explicit} versus when it is computed in a centralized manner, i.e., when solving the same optimization problem with the optimization solver CVX \cite{grant_graph_2008,grant_cvx_2013}, in order to illustrate that the performance is optimal.

\begin{figure}[htp]
    \centering
    \subcaptionbox{Runtime for each of the three cases for varying system sizes.}{ 
\includegraphics[width=0.7\columnwidth]{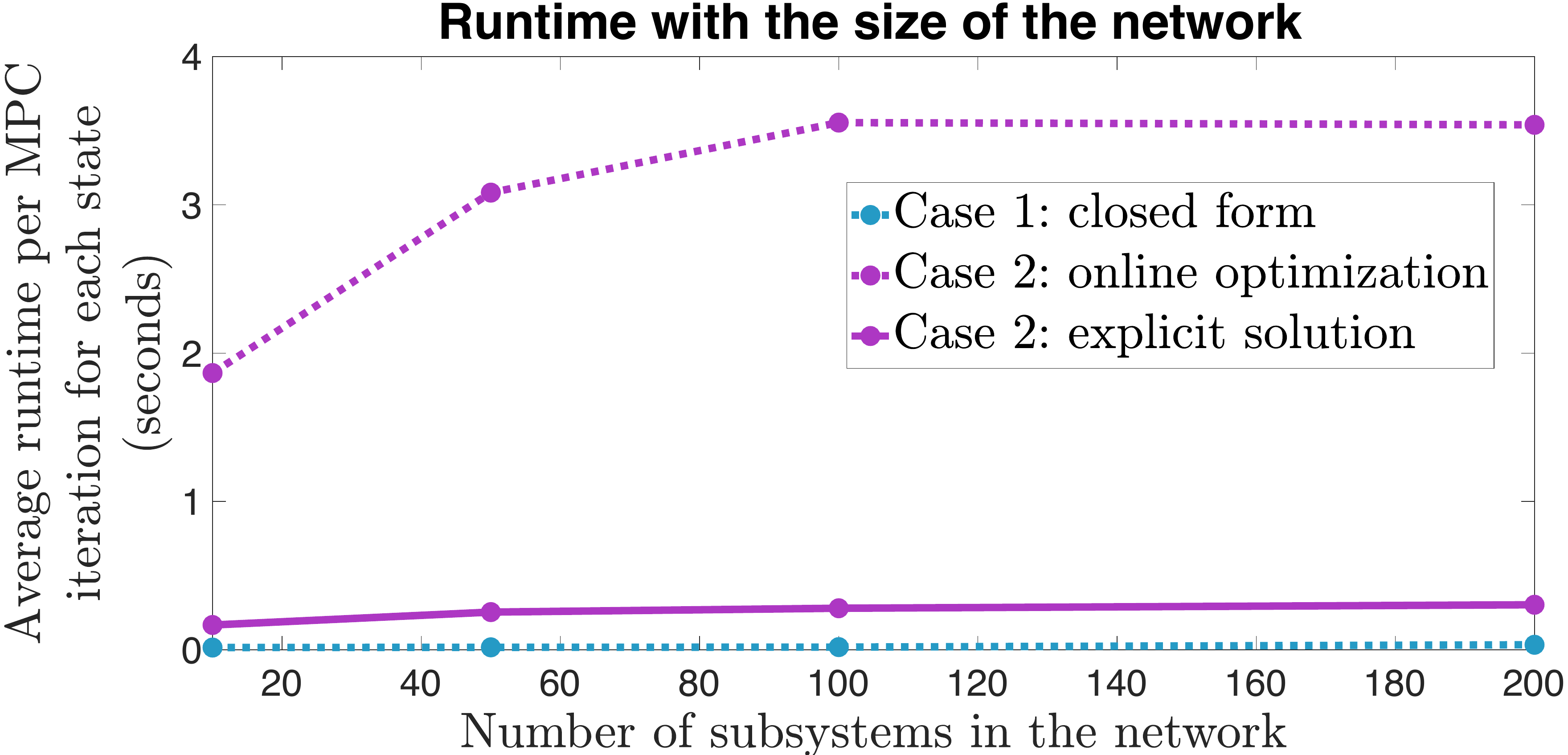}
}\vfill
    \subcaptionbox{Comparison of the optimal cost when solved with Algorithm \ref{alg: I explicit} versus the centralized problem with CVX for different network sizes.}{
\includegraphics[width=0.7\columnwidth]{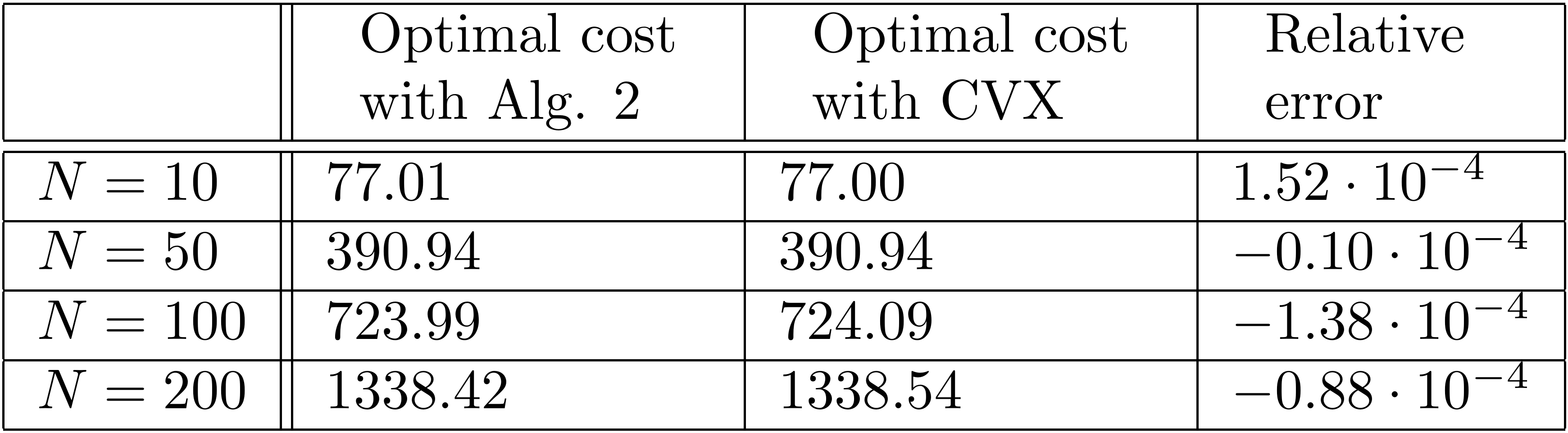}
}\vfill
    \caption{Simulation results.}
    \label{fig: scalability}
\end{figure}

These simulations illustrate that the proposed explicit MPC scheme scales gracefully with the size of the system, can be applied to arbitrarily large systems as long as they satisfy the required structural assumptions, and achieves optimal solutions. 
\section{Conclusion}\label{sec:conclusion}
We introduced an explicit solution to the MPC problem that can be applied to large networked systems.  Together with some natural separability assumptions on the objective function and constraints, leads to an \emph{explicit} distributed and localized synthesis and implementation of the MPC controller derived in~\cite{amo_alonso_distributed_2019}.  

Inspired by \cite{bemporad_explicit_2000}, our explicit solution partitions the space into three regions per state/input instantiation, and with the assumptions  that no coupling between states is allowed, each subsystem can solve for each state (input) instantiation sequentially, which results in a fast computation runtime per subsystem. Since each subsystem runs its own optimization problem in parallel, this results in large runtime improvements. Computational experiments show that the runtime of each MPC iteration per subsystem in the network scales with $O(1)$ complexity as the size of the network increases.

\section*{Appendix A: System Level Synthesis }\label{Appendix A}

The main result of the SLS theory is stated below. The reader is referred to Theorem 2.1 in \cite{anderson_system_2019} for a detailed explanation and proof. 
\begin{theorem}{\label{thm: SLS}}
For the dynamical system~\eqref{eq:SLScompact} evolving over a finite horizon, under the state-feedback policy $\mathbf u = \mathbf K \mathbf x$, for $\mathbf{K}$ a block-lower-triangular matrix, the following are true
\begin{enumerate}
    \item The affine subspace
    \begin{equation}\label{eq: constraint}
        \left[I-Z\mathcal A\ \ -Z\mathcal B\right]\left[\begin{array}{c}\mathbf{\Phi}_{x}\\\mathbf{\Phi}_{u}\end{array}\right] = I
    \end{equation}
    parameterizes all possible system responses \eqref{eq:response}.
    
    \item For any block-lower-triangular matrices $\left\{\mathbf{\Phi}_{x},\mathbf{\Phi}_{u}\right\}$ satisfying (\ref{eq: constraint}), the controller $\mathbf{K} = \mathbf{\Phi}_{u}\mathbf{\Phi}_{x}^{-1}$ achieves the desired response \eqref{eq:response} from $\mathbf w \mapsto (\mathbf x,\mathbf u)$.
\end{enumerate}
\end{theorem}
As $\mathbf{\Phi}_{x}$ and $\mathbf{\Phi}_{u}$ appear in an affine manner in~\eqref{eq: constraint} they can be incorporated into a convex program.

\section*{Appendix B: Structured system example} \label{Appendix B}
\begin{figure}[t]
    \centering
    \includegraphics[width=0.70\columnwidth]{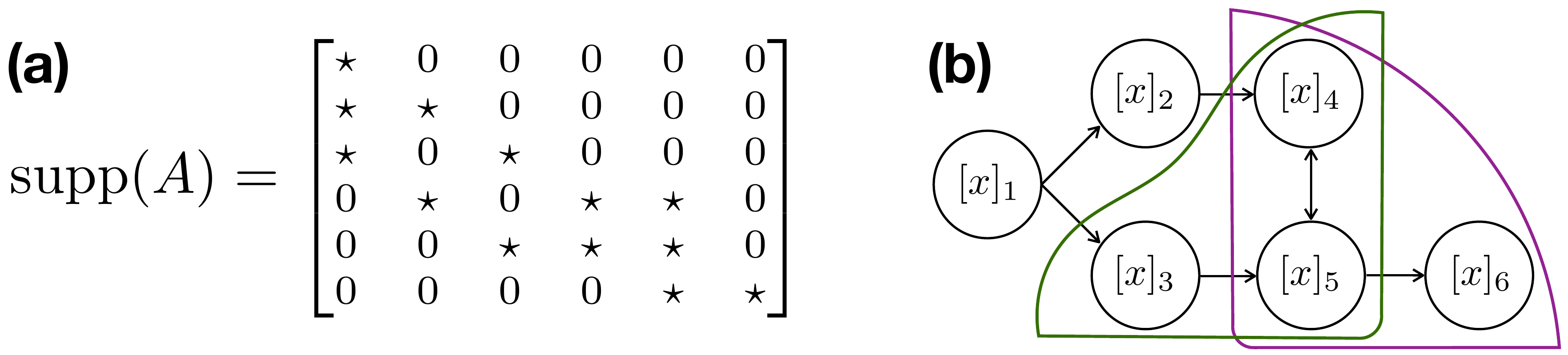}
    \caption{{(a) Support of matrix $A$. (b) Example of $1$-incoming and $1$-outgoing sets for subsystem $5$.}}
    \label{fig: in_out}
\end{figure}

This example was directly taken from \cite{amo_alonso_distributed_2019}. Consider a system \eqref{eq: LTI system} composed of $N=6$ scalar subsystems, with $B = I_6$ and $A$ matrix with support represented in Figure \ref{fig: in_out}(a).
This induces the interconnection topology graph $\mathcal G_{(A,B)}$ illustrated in Figure \ref{fig: in_out}(b). As shown, the $d$-incoming and $d$-outgoing sets can be directly read off from the interaction topology.  For example, for $d=1$, the $1$-hop incoming neighbors for subsystem $5$ are subsystems $3$ and $4$, hence $\textbf{in}_{5}(1)=\{3,4\}$; similarly, we observe that $\textbf{out}_{5}(1)=\{4,6\}$.

\bibliographystyle{IEEEtran}
\vspace{-5mm}
\bibliography{SLS}

\end{document}